\documentclass[reqno,10pt]{amsart}
\usepackage{amsmath}
\usepackage{amsthm}
\usepackage{amssymb}
\usepackage{bm}
\usepackage{pstricks}

\newcommand{\fb}{{\mathbf f}}

\newcommand{\LL}{\mathcal L}
\newcommand{\nb}{{\mathbf n}}

\newcommand{\qb}{{\mathbf q}}

\newcommand{\Qb}{{\mathbf Q}}

\newcommand{\RR}{\mathbb R}
\renewcommand{\SS}{{\mathcal S}}

\newcommand{\vb}{{\mathbf v}}

\newcommand{\Vb}{{\mathbf V}}

\newcommand{\ZZ}{{\mathbb Z}}

\newcommand{\str}{\widehat{\eps}}
\newcommand{\rstr}{\widetilde{\eps}}

\newcommand{\eps}{\varepsilon}
\newcommand{\SLZ}[1]{\text{SL}(#1;\ZZ)}

\DeclareMathOperator{\diag}{diag}

\theoremstyle{plain}

\newtheorem{assumption}{Assumption}[section]

\theoremstyle{remark}
\newtheorem{remark}{Remark}[section]

\title[Periodic boundary conditions for nonequilibrium flows]{Periodic boundary conditions for long-time nonequilibrium molecular
dynamics simulations of incompressible flows}
\author{Matthew Dobson}\address{Matthew Dobson, Department of Mathematics and Statistics,
710 N.~Pleasant Street,
University of Massachusetts,
Amherst, MA 01003-9305, USA
}
\email{dobson@math.umass.edu}

\date{\today}
\begin{document}
\maketitle

\begin{abstract}
This work presents a generalization of the Kraynik-Reinelt (KR) boundary
conditions for nonequilibrium molecular dynamics simulations.  In the
simulation of steady, homogeneous flows with periodic boundary conditions, the
simulation box moves with the flow, and it is possible for particle replicas to
become arbitrarily close, causing a breakdown in the simulation.  The KR
boundary conditions avoid this problem for planar elongational flow and general
planar mixed flow [J. Chem. Phys {\bf 133}, 14116 (2010)] through careful choice
of the initial simulation box and by periodically remapping the simulation box
in a way that conserves replica locations.  In this work, the ideas are
extended to a large class of three dimensional flows by using multiple
remappings for the simulation box.  The simulation box geometry is no longer
time-periodic (which was shown to be impossible for uniaxial and biaxial
stretching flows in the original work by Kraynik and Reinelt [Int. J.
Multiphase Flow {\bf 18}, 1045 (1992)]).  The presented algorithm applies to all
flows with nondefective flow matrices, and in particular, to uniaxial and
biaxial flows.
\end{abstract}

\section{Introduction}

Nonequilibrium molecular dynamics techniques are widely employed
in the study of molecular fluids under steady flow.  Periodic boundary
conditions (PBCs) are employed to study bulk properties of a fluid, but
standard PBCs with a fixed simulation box are incompatible with a homogeneous
linear background flow $A = \nabla u \in \RR^{3 \times 3}$, such as shear or
elongational flow.   In such a simulation, the periodic replicas of a particle
have different velocities, consistent with the background flow.  If we let 
$$
L_t = \bigg[ \vb^1_t \ \vb^2_t \ \vb^3_t \bigg] \in \RR^{3 \times 3},
t \in [0, \infty)
$$ 
denote the time-dependent lattice basis vectors defining the simulation box, 
then a particle with phase coordinates $(\Qb, \Vb)$ has periodic replicas
with coordinates at $(\Qb + L_t \nb, \Vb + A L_t \nb)$ for
all integer triples $\nb \in \ZZ^3.$  The velocity relations 
\begin{equation*}
\frac{d}{dt} (\Qb + L_t \nb) = \Vb + A L_t \nb \text{ for all } \nb \in \ZZ^3
\end{equation*} imply that
the simulation box must move with the flow, 
\begin{equation}
\label{Lt}
\frac{d}{dt} L_t = A L_t, \text{ which has solution } L_t = e^{A t} L_0.
\end{equation}   
For general flows, depending on the orientation of $L_0$ the simulation box can
become quite elongated so that a particle is approached by its periodic
replicas, which causes numerical instability in the simulation.   For example,
a planar elongational flow whose contraction is parallel to one of the
simulation box edges $\vb^i_0$ has one periodic direction that shrinks
exponentially fast.  This puts a finite limit on the simulation
stability~\cite{houn92, bara95}.  While these time periods are sometimes long
enough to allow for the accurate computation of statistical observables in
simple molecular fluids, there is need for boundary conditions without time
limitations for the simulation of complex molecular systems. 
  
For shear flow, the Lees-Edwards boundary conditions~\cite{lees72} allow for
time-periodicity in the deforming simulation box itself.  For planar
elongational flow, the Kraynik-Reinelt (KR) boundary conditions~\cite{kray92,
todd98, todd99, bara99} achieve time periodicity in the simulation box by
carefully choosing the vectors defining the initial simulation box.  In
particular, the box is rotated so that the edges form an angle of approximately
31.7 degrees with respect to the background flow.  However, the KR formalism
does not apply to general three dimensional flows, in particular it cannot
treat uniaxial or biaxial flow~\cite{kray92}.  In this paper, we generalize the
KR boundary conditions to handle any homogeneous, incompressible,
three-dimensional flow whose velocity gradient is a nondefective matrix (see
Section~\ref{sec:model} for a precise description of the flow types handled).
We greatly enlarge the class of flows handled, including uniaxial and biaxial
flows.  The proposed algorithm gives an initial orientation for the lattice
vectors $L_0,$ evolves the vectors according to the differential
equation~\eqref{Lt}, and remaps the vectors in a fashion that preserves the
periodic lattice structure and keeps the total deformation bounded for all
time.  Unlike Lees-Edwards and Kraynik-Reinelt boundary conditions, the
boundary conditions do not in general have a time-periodic simulation box;
however, the deformation of the simulation box is kept bounded and particle
replicas stay separated by a bounded distance.  

In Section~\ref{sec:kr} we review the KR boundary conditions and describe them
in a framework useful for the generalization later.  In Section~\ref{sec:genkr}
the new boundary conditions are derived and explained theoretically.
Section~\ref{sec:algo} contains a self-contained description of the algorithm
with default choices for parameters given.

We note that the boundary conditions described here are not tied to a
particular choice of nonequilibrium dynamics.  Typically, the flow in a
nonequilibrium simulation is driven by a specialized dynamics, for example, the
deterministic SLLOD~\cite{evan07,edbe86} or g-SLLOD~\cite{tuck97,edwa06}
dynamics or the nonequilibrium stochastic dynamics such as those
in~\cite{mcph01,dobs12}.  

\section{Flow Types and Automorphisms}
\label{sec:model}

Since the background flow treated here is incompressible, $A$ is a trace-free
matrix.  Let $J = S^{-1} A S$ denote the real Jordan canonical form for $A,$
where all 3 by 3 matrices fall in four possible cases, 
\begin{align}
\label{eq:J_nondef}J_1 = \left[ 
\begin{array}{rrr}
\eps_1 & 0 & 0 \\
0 & \eps_2 & 0 \\
0 & 0 & - \eps_1 - \eps_2 
\end{array}\right], 
J_2 = \left[ \begin{array}{rrr}
\eps & -r & 0 \\
r & \eps & 0 \\
0 & 0 & -2 \eps 
\end{array}\right], \\
\label{eq:J_def}
J_3 = \left[ 
\begin{array}{rrr}
\eps & 1 & 0 \\
0 & \eps & 0 \\
0 & 0 & - 2 \eps 
\end{array}\right], \text{ or }
J_4 = \left[ \begin{array}{rrr}
0 & 1 & 0 \\
0 & 0 & 1 \\
0 & 0 & 0 
\end{array}\right].
\end{align}
The form $J_1$ includes several standard matrices, for example planar
elongational flow (PEF) where $\eps_1 = -\eps_2,$ uniaxial stretching flow
(USF) where $\eps_1 = \eps_2 < 0,$ and biaxial stretching flow (BSF) where
$\eps_1 = \eps_2 > 0.$  The matrix $J_2$ arises in the case of complex
eigenvalues, corresponding to a rotational flow (which may be an inward spiral
$a < 0,$ an outward spiral $a > 0,$ or a center $a=0$).  Both $J_3$ and $J_4$
are defective matrices, since they have rank-deficient eigenspaces.  The
generalized KR boundary conditions apply to any matrices of the form $J_1$ or
$J_2$.  For $J_3,$ if $\eps = 0,$ then this is a case of planar shear flow and
the Lees-Edwards boundary conditions can be employed.  Likewise, similar
boundary conditions can be employed for the case $J_4;$ however, we have not
been able to extend the boundary conditions described here to the $J_3$ case
for nonzero $\eps.$

In the following, we transform the lattice $L_t$ with elements of $\SLZ{3},$
the matrix group of orientation-preserving linear lattice automorphisms.  This
is the set of all three by three matrices with integer entries whose
determinant is 1.  By Cramer's rule, such a matrix has an inverse with integer
entries.  For any $M \in \SLZ{3},$ the lattices generated by $L_t$ and $L_t M$
are identical, and thus, the two sets of particles $\{ \Qb_i + L_t \nb \, | \,
\nb \in \ZZ^3\}$ and $\{ \Qb_i + L_t M \nb \, | \, \nb \in \ZZ^3 \}$ are
identical.  Applying such an automorphism transforms the simulation box without
changing the simulated dynamics.  Through the careful choice of initial
simulation box $L_0$ and automorphisms, we can simulate a system where all
particles maintain a minimum distance from their periodic replicas for all
time.

\section{KR boundary conditions and planar flows}
\label{sec:kr}

We first present a review of the KR boundary conditions for planar
elongational flow along with a summary of techniques for
other planar flows.

\subsection{KR boundary conditions for planar elongational flow}
Consider a diagonal flow of the form
\begin{equation*}
A = \left[ \begin{array}{rrr} 
{\eps}&0 &0\\
0&-{\eps} &0\\
0&0&0
\end{array}\right],
\end{equation*} 
where ${\eps} >0$.  

The KR boundary conditions~\cite{kray92} consist in choosing a basis
for the unit cell such that after a finite time, the elongational flow
maps the lattice generated by the unit cell onto itself.  That is, one
finds a basis $L$ and time $t_* > 0$ such that $$e^{A t_*} L = L M,$$ for
some  $M \in \SLZ{3}.$  The mapping $M$ is a parameter of the
algorithm. 
The method was first described in~\cite{kray92}, where the authors
showed how to find reproducible square and hexagonal lattices in
planar elongational flow. 
In~\cite{todd98,bara99,todd99} the authors employed these reproducible
lattices in nonequilibrium molecular dynamics simulations by using
them to describe the periodicity of groups of particles.

Choose $M \in \SLZ{3}$ with positive eigenvalues, other than the identity
matrix.  
For example, the
choice 
$$
M = \left[
\begin{array}{rrr}
 2 & -1 &  0 \\
-1 &  1 &  0 \\
 0 &  0 &  1
\end{array}
\right]$$
is common, and it
has been shown to give a system with the largest possible minimal
spacing between periodic replicas~\cite{kray92}.
Let
$V$ denote a matrix of eigenvectors for $M$, and let
$$\Lambda = \left[ 
\begin{array}{rrr}
\lambda & 0 & 0 \\
0 & \lambda^{-1} & 0 \\
0 &    0         & 1 
\end{array}
\right]
$$ denote the matrix of
corresponding eigenvalues, so that
\begin{equation*}
M V = V \Lambda.
\end{equation*} 
We order the eigenvalues so that $\lambda > 1.$
The fact that the eigenvalues are inverses of one another follows from
$\det(M) = 1.$ 
We define the lattice time period  
\begin{equation} 
\label{tzero}
t_* = \frac{\log(\lambda)}{{\eps}}
\end{equation}  
so that $e^{{\eps} t_*} =
\lambda.$
Let $L_0 = V^{-1}$ be the matrix of initial lattice vectors.
Note that while it is typical to choose eigenvectors to have norm one, 
the vectors in $V$ should be scaled so
that the unit cell $L_0 = V^{-1}$ has the desired volume for the simulation box.
If one chooses the vectors of $V$ to have the same length, then the
vectors of $V$ are orthogonal, and $V^{-1} = \frac{1}{\det(V)^2} V^T.$
Since the lattice vectors move with the flow as in~\eqref{Lt}, at time
$t_*$, they satisfy 
\begin{equation*}
L_{t_*} = e^{A t_*} V^{-1} = \Lambda V^{-1}
= L_0 M.
\end{equation*}  
Thus, the lattice vectors $L_{t_*}$ generate the same lattice
as $L_0$, demonstrating the time periodicity of the lattice.  In
simulations, the simulation box is remapped by setting  
$$L_{t_*^+} := L_0$$ 
to
avoid the use of highly elongated basis vectors.  This transformation
does not move any of the periodic replicas of the particles in the
simulation; however, since the basis vectors have changed, the
periodic boundary conditions need to be applied on stored particle positions
so that the stored
particle displacements fall within the simulation box.  

\subsection{General planar flows}
As mentioned in~\cite{kray92} and implemented for mixed flow in~\cite{hunt10},
the above algorithm can be applied to certain nondiagonal matrices $A$ by
diagonalization.  However, in~\cite{kray92}, it is shown by consideration of
the characteristic polynomial for members of $\SLZ{3}$ that there is no
reproducible lattice for either USF or BSF.  Suppose now that $A$ denotes a
general incompressible planar flow, that is, all nonzero entries of the matrix
act on a two-dimensional eigenspace.  This corresponds to cases $J_1$ with
$\eps_1 = \eps_2$, $J_2$ with $\eps = 0,$ or $J_3$ with $\eps=0$
in~\eqref{eq:J_nondef} and~\eqref{eq:J_def}.  There are three cases to
consider, two nonzero real eigenvalues, two purely imaginary eigenvalues, or
only zero eigenvalues.  

\subsubsection{Elongational flow}
\label{sec:2D_elong}
If the eigenvalues of $A$ are real and distinct then 
$A$ is diagonalizable and corresponds to an elongational
flow.  Let $S$ denote a matrix of eigenvectors and $D$ denote the
matrix of eigenvalues for $A$ so that $A S = S D.$ Then, upon choosing
the basis $L_0 = S V^{-1},$ we have 
\begin{equation*}
\begin{split} 
L_{t_*} &= e^{A t_*} L_0
= e^{A t_*} S V^{-1}
= S e^{D t_*} V^{-1} 
= S V^{-1} M 
= L_0 M.
\end{split} 
\end{equation*}
We note that this includes the mixed flow case treated in~\cite{hunt10}.

\subsubsection{Rotational flow}
In the case the eigenvalues are pure imaginary, and the
flow is rotational.  Writing $A$ in real Jordan normal
form, we choose real $S$ so that 
\begin{equation*}
S^{-1} A S = \left[\begin{array}{rrr}
 0 & r & 0\\
-r & 0 & 0\\
 0 & 0 & 0
\end{array}\right].
\end{equation*} 
We define $L = S$ and then have that
$e^{A t} L =  R_t L,$ where $R_t$ is a rotation for all $t$.  
There is no need to reset the simulation box in this case.

\subsubsection{Shear flow}
\label{sec:2D_shear}
The final case of all zero eigenvalues corresponds to shear flow.  
We note that in this case, there is a $t_*$
and $S$ such that
$e^{A t_*} S = S M,$
for \begin{equation*}
M = \left[ \begin{array}{rrr}
1 & 1 & 0\\
0 & 1 & 0\\
0 & 0 & 1
\end{array} \right].
\end{equation*}
This is the Lagrangian rhomboid scheme, which is equivalent to the
Lees-Edwards boundary conditions~\cite{lees72, evan07}.

\section{Generalized KR boundary conditions}
\label{sec:genkr}
In this section we generalize the boundary conditions to nondefective
incompressible linear flows in three dimensions.  In the following,
rather than find a time $t_0$ such that $L_{t_0} = L_0 M$ for a single
automorphism $M \in \SLZ{3},$ we consider the successive application
of two different automorphisms $M_1, M_2 \in \SLZ{3}$ to $L_t$ in order to
keep the total deformation of the unit cell small for all times. 

Suppose that  $M_1, M_2 \in \SLZ{3}$ are a pair of commuting,
symmetric automorphisms.  Then the matrices are simultaneously
diagonalizable by an orthogonal matrix $V.$  Let 
\begin{equation*} 
\Lambda_{i}=V^{-1} M_i V
\end{equation*} 
denote the matrix of eigenvalues corresponding to $M_i,$ whose
diagonal entries are denoted by $ \lambda_{i, 1}, \lambda_{i, 2},
\lambda_{i, 3}.$  We define the logarithm of the ordered spectrum for
each operator
\begin{equation}
\label{logspec}
\hat{\omega}_i = \left[ 
\begin{array}{c}
\log \lambda_{i,1} \\
\log \lambda_{i,2} \\
\log \lambda_{i,3} \\
\end{array}
\right]
= \log \diag( V^{-1} M_i V ),
\end{equation} 
where $\diag(M)$ denotes the column vector made up of the diagonal
entries of the matrix $M.$  We assume the following about $M_1$ and
$M_2.$
\begin{assumption}
\label{spec_ass}
We assume that $M_1, M_2 \in \SLZ{3}$ are symmetric, commute, and have
positive eigenvalues.  We assume that $\hat{\omega}_1$ and
$\hat{\omega}_2,$ defined in~\eqref{logspec}, are linearly
independent.
\end{assumption}

An example of such a pair of matrices is given in
Section~\ref{sec:algo}.  Note that the choice of $M_1$ and $M_2$ does
not depend on the matrix $A$.  

We describe the technique first in the
diagonal case before discussing in turn the four possible cases for
three dimensional flows.  After the derivation given here, the main
algorithm is presented in a concise form in Section~\ref{sec:algo}.

\subsection{Diagonal case}
Let us first consider a diagonal flow of the form
\begin{equation}
\label{adiag}
A = \left[ \begin{array}{rrr}  
\eps_1 & & \\
& \eps_2 & \\
& & \eps_3 \\
\end{array}
\right],
\end{equation}
where $\eps_1 + \eps_2 + \eps_3 = 0.$  Then the matrix
exponential
\begin{equation}
\label{expAt}
e^{A t} = \left[ \begin{array}{rrr}  
e^{\eps_1 t} & & \\
& e^{\eps_2 t} & \\
& & e^{\eps_3 t} \\
\end{array}
\right],
\end{equation}
is diagonal for all time $t.$

Let $M_1$ and $M_2$ satisfy Assumption~\ref{spec_ass}.  We choose
initial lattice basis $L_0 = V^{-1},$ where $V$ diagonalizes $M_1$ and
$M_2.$   Applying the transformation $M_i$ to $L_t$ gives  
\begin{equation*}
\begin{split} 
L_t M_i &= e^{A t} V^{-1} M_i \\
&= \left[ 
\begin{array}{rrr} 
e^{\eps_1 t} & & \\
& e^{\eps_2 t} & \\
& & e^{\eps_3 t}
\end{array}
\right] 
\left[ 
\begin{array}{rrr} 
\lambda_{i, 1} & & \\
& \lambda_{i, 2} & \\
& & \lambda_{i, 3}
\end{array}
\right] V^{-1} \\
&= \exp\left( \left[ 
\begin{array}{rrr} 
{\eps_1 t} + \log\lambda_{i, 1} & & \\
& {\eps_2 t} + \log\lambda_{i, 2} & \\
& & {\eps_3 t} + \log\lambda_{i, 3} 
\end{array}
\right] \right)  V^{-1}.
\end{split} 
\end{equation*} 
Similarly, if we apply multiple transformations at once, we have
\begin{equation} 
\label{mult_trans}
\begin{split} 
L_t M_1^{n_1} M_2^{n_2} 
= \exp\left( 
\left[ 
\begin{array}{rrr} 
{\eps_1 t} & & \\
& {\eps_2 t} & \\
& & {\eps_3 t} 
\end{array}
\right]
+ \sum_{i=1}^2 n_i 
\left[ 
\begin{array}{rrr} 
\log\lambda_{i, 1} & & \\
& \log\lambda_{i, 2} & \\
& & \log\lambda_{i, 3} 
\end{array}
\right] \right)  V^{-1}
\end{split}
\end{equation} 
where $n_1, n_2 \in \ZZ.$  The idea of the algorithm presented in
Section~\ref{sec:algo} is to apply automorphisms so that the
argument of the exponential in~\eqref{mult_trans} 
stays  bounded for all times $t >0.$

We define a vector that equals the diagonal part of the stretch, 
\begin{equation*} 
\str_t = \left[ 
\begin{array}{c} 
\eps_1 t \\ 
\eps_2 t \\
\eps_3 t 
\end{array} \right],
\end{equation*} 
and note that $\str_t,$ $\hat{\omega}_1,$ and $\hat{\omega}_2$ belong
to the two dimensional subspace $\SS \subset \RR^3$ of mean-zero
vectors.  The vectors $\hat{\omega}_1$ and $\hat{\omega}_2$ generate a
lattice in $\SS,$ 
\begin{equation*}
\LL = \left\{ \left(n_1 - \frac{1}{2}\right) \hat{\omega}_1 + \left(n_2
- \frac{1}{2}\right) \hat{\omega}_2 \ | \ 
n_1, n_2 \in \ZZ \right\},
\end{equation*}  
where we have added an offset of $1/2$ so that the unit cell
$$
\widehat{\Omega} = \left\{ \theta_1 \hat{\omega}_1 + \theta_2 \hat{\omega_2} \ | \ 
\theta_1, \theta_2 \in \left(-\frac{1}{2},\frac{1}{2}\right] \right\}
$$ 
is centered at the origin.  At each time $t > 0,$ by applying powers of the
automorphisms to the lattice, we can transform so that the remapped simulation
box
$$\widetilde{L}_t = L_t M_1^{n_1} M_2^{n_2}$$ 
has a small stretch vector
$\tilde{\eps}_t = \hat{\eps}_t + n_1 \hat{\omega_1} + n_2 \hat{\omega_2}.$

\subsection{Diagonalizable flow}

Suppose that $A$ is diagonalizable, 
\begin{equation*}
A = S D S^{-1}.
\end{equation*} 
As pointed out for the planar case in Section~\ref{sec:2D_elong}, we can extend
the above algorithm, by choosing $L_0 = S V^{-1}.$ We then have
\begin{equation*}
L_t M_1^{n_1} M_2^{n_2} = e^{A t} S V^{-1} M_1^{n_1} M_2^{n_2}
= S e^{D t} V^{-1} M_1^{n_1} M_2^{n_2}.
\end{equation*} 
The automorphisms act to bound the stretch vector corresponding to the diagonal term
$e^{D t}.$ 
We note that since $S$ is not orthogonal if $A$ is nonsymmetric,  the original
lattice vectors $L_0$ are not orthogonal in that case.

\subsection{Complex eigenvalues}

It is also possible that $A$ has a pair of complex eigenvalues and a single
real eigenvalue.  We denote the spectrum of A as $\{ \eps + i r, \eps - i r, -2
\eps\}.$  In this case, we write the real Jordan normal form for the matrix,
\begin{equation*}
A = S J_2 S^{-1},
\end{equation*} 
where $S$ is real and $J_2$ is the block-diagonal matrix  
\begin{equation*}
J_2 = \left[
\begin{array}{rrr}
\eps & r & 0 \\
-r & \eps & 0 \\
0 & 0 & -2\eps 
\end{array}
\right].
\end{equation*} 
We decompose $J_2 = D + B$
where 
\begin{equation*}
D = \left[
\begin{array}{rrr}
\eps & 0 & 0 \\
0 & \eps & 0 \\
0 & 0 & -2\eps 
\end{array}
\right] \text{ and } 
B = \left[
\begin{array}{rrr}
0 & r & 0 \\
-r & 0 & 0 \\
0 & 0 & 0 
\end{array}
\right].
\end{equation*} 
We note that since $D B = B D,$ the matrix exponential splits into a rotation
and a stretch giving 
\begin{equation*} 
e^{A t} = S e^{J t} S^{-1} = S e^{B t} e^{D t} S^{-1}, 
\end{equation*}
where $e^{B t}$ is a rotation matrix. We again take initial lattice vectors 
$L_0 = S V^{-1}$ and control size of the stretch vector  
\begin{equation*}
\str_t = 
\left[ \begin{array}{r} 
\eps t\\
\eps t\\
-2 \eps t 
\end{array} \right],
\end{equation*} 
using the automorphisms $M_1$ and $M_2.$ No effort is made to undo the effect
of $e^{B t}$ since it is simply a rotation.  

\subsection{Defective matrices}

The final possible case is when $A$ is a defective matrix, that is, it has a
repeated eigenvalue whose eigenspace does not have full rank.  In three
dimensions, a defective matrix can only occur for a matrix with a real
spectrum, and so the only possible Jordan forms, up to rearrangement of the
blocks, are
\begin{equation*}
J_3 = \left[ 
\begin{array}{rrr}
\eps & 1 & 0 \\
0 & \eps & 0 \\
0 & 0 & -2\eps 
\end{array}\right] \text{ or }
J_4 = \left[ \begin{array}{rrr}
0 & 1 & 0 \\
0 & 0 & 1 \\
0 & 0 & 0 
\end{array}\right].
\end{equation*}

We can treat the $J_4$ case very similarly to the shear flow case in
Section~\ref{sec:2D_shear}, using the identity
\begin{equation*}
e^{J_2 t} = \left[ \begin{array}{rrr}
1 & t & \frac{t^2}{2} \\
0 & 1 & t \\
0 & 0 & 1 
\end{array}\right].
\end{equation*}   
We choose initial lattice basis $L_0 = S$ and note that at time $t_0 = 2,$ we
have
\begin{equation*}
\begin{split}
L_{t_0} &= e^{2 A} S \\
&= S e^{2 J_2} \\
&= S \left[
\begin{array}{rrr}
1 & 2 & 2 \\
0 & 1 & 2 \\
0 & 0 & 1
\end{array}
\right]\\
&= S M,
\end{split} 
\end{equation*} 
where $M \in \SLZ{3}.$

 We have not been able to generalize our algorithm to
the case of $J_3,$ when $\eps \neq 0.$  The difficulty lies with the off-diagonal
terms of the matrix exponential
$$
e^{ J_3 t} = \left[
\begin{array}{rrr}
e^{\eps t} & t e^{\eps t} & 0 \\
0 & e^{\eps t} & 0 \\
0 & 0 & e^{-2 \eps t}
\end{array}
\right].
$$
One approach considered is to find 
matrices $M_j \in \SLZ{3}$ and a common matrix $V$ such that $V^{-1} M_j V$ is upper triangular,
in order to control the diagonal and off-diagonal terms at the same time,
but we have not had success in such a construction.

\section{Algorithm}
\label{sec:algo}
We now provide an explicit construction of the generalized KR boundary
conditions algorithm.  The following two matrices are in $\SLZ{3}$ and they
commute:
\begin{equation*}
M_1 = \left[ 
\begin{array}{rrr} 
1 & 1 & 1 \\
1 & 2 & 2 \\
1 & 2 & 3 
\end{array}
\right] \qquad 
M_2 = \left[ 
\begin{array}{rrr} 
 2 & -2 &  1 \\
-2 &  3 & -1 \\
 1 & -1 &  1 
\end{array}
\right]. 
\end{equation*}
We choose the initial lattice vectors $L_0 = a V^{-1},$ where $V$ denotes the
matrix of eigenvectors for $M_1$ and $M_2,$ and $a^3$ is the volume of the
simulation box.  We fix the choice of ordering for the eigenvectors by giving
the first few digits of $V^{-1}$, 
\begin{equation*}
V^{-1} = \left[  
\begin{array}{rrr} 
0.591 & -0.737 & 0.328 \\ 
0.737 &  0.328 &-0.591 \\ 
0.328 &  0.591 & 0.737 
\end{array}
\right]
\end{equation*}
Direct computation shows that the ordered spectra of the two operators
are positive and the corresponding $\hat{\omega}_i,$ given by
\begin{equation*}
\hat{\omega}_1 \approx \left[
\begin{array}{r}
  - 1.178  \\
    1.619  \\
  - 0.441  
\end{array}
\right]
\qquad
\hat{\omega}_2 \approx \left[
\begin{array}{r}
    1.619  \\
  - 0.441  \\
  - 1.178  \\
\end{array}
\right]
\end{equation*} 
are linearly
independent. 

Suppose that $A$ is written in real Jordan normal form $A = S J S^{-1}$ and $J$
is decomposed as $J = D + B$ where 
\begin{equation*}
D = \left[
\begin{array}{rrr}
\eps_1 & 0 & 0 \\
0 & \eps_1 & 0 \\
0 & 0 & \eps_3 
\end{array}
\right] \text{ and } 
B = \left[
\begin{array}{rrr}
0 & r & 0 \\
-r & 0 & 0 \\
0 & 0 & 0 
\end{array}
\right].
\end{equation*}   
This encompasses both diagonalizable flow (where $B = 0$) and the case of
complex eigenvalues, but does not include the defective matrix
case~\eqref{eq:J_def}. 
 
\begin{figure}[tb]
\centerline{\input{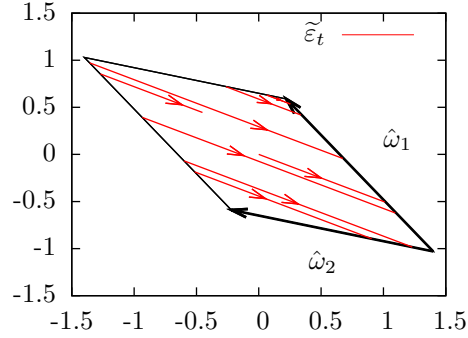}}
\caption{\label{fig:stretch}  As the simulation progresses, $\rstr_t$ 
traces a curve in the unit cell $\hat{\Omega}$ within $\SS.$
Here, the unit cell $\widehat{\Omega}$ of the lattice in stretch space 
has been projected 
into the xy plane.  The lines within the parallelogram 
denote the evolution of $\rstr_t$ during
a simulation of uniaxial stretching flow.  The depicted unit cell corresponds
to the example automorphisms given in Section~\ref{sec:algo}.}
\end{figure}

For time $t \geq 0$, we define the reduced stretch $\rstr_t$ as follows
\begin{equation*}
\frac{d}{d t} \, \rstr_t = \left[
\begin{array}{r}
\eps_1 \\
\eps_2 \\
\eps_3 
\end{array}
\right], \qquad
\rstr_0 = \left[
\begin{array}{r}
0 \\
0 \\
0
\end{array}
\right],
\end{equation*}
where $\rstr_t$ is restricted to be within the unit cell 
\begin{equation*}
\widehat{\Omega} = \left\{ \theta_1 \hat{\omega}_1 + \theta_2
\hat{\omega}_2 \  | \  
 \theta_1, \theta_2 \in \left(-\frac{1}{2}, \frac{1}{2}\right] \right\},
\end{equation*} 
by periodic boundary conditions.  An example curve $\rstr_t$ is depicted in
Figure~\ref{fig:stretch}.  The lattice basis vectors for the simulation are
then defined to be 
\begin{equation*}
\widetilde{L}_t = S e^{B t} e^{\rstr_t} V^{-1},
\end{equation*} 
where we define 
$$e^{\rstr_t} = \exp\left( \left[ \begin{array}{rrr}
\rstr_{t,1} & & \\
& \rstr_{t,2} & \\
& & \rstr_{t,3} 
\end{array} \right] \right).$$
This process can be carried out for arbitrarily long times, and the stretch
$\rstr_t$ stays bounded for all times.

This gives the following pseudocode for the discretized version
of the NEMD system:  \\ 
Given $S, D, B,$ and the time step $\Delta t$,  compute $(\delta_1, \delta_2)$ so that $\delta_1 \hat{\omega}_1 +
\delta_2 \hat{\omega}_2 = [\eps_1, \eps_2, \eps_3]^T.$  \\
For each time step do:
\begin{enumerate}
\item $\theta_i \leftarrow \theta_i + \delta_i \Delta t$
\item $\theta_i \leftarrow \theta_i - \mathrm{round}(\theta_i)$
\item $\tilde{\eps}_t \leftarrow \theta_1 \hat{\omega_1} + \theta_2 \hat{\omega_2}$
\item $\widetilde{L}_t \leftarrow S e^{B t} e^{\rstr_t} V^{-1}$
\end{enumerate}

Note that we recompute the lattice basis vectors at each step,
and we do not explicitly apply automorphisms nor do we directly
reset the lattice vectors.

\subsection{Minimum replica distance}

The boundary conditions above limit the stretch $\rstr_t$ to
live within the unit cell $\widehat{\Omega}$ which is defined by the
 vectors~\eqref{logspec}. 
The minimum distance
between a particle and a periodic replica within the simulation is
given by
\begin{equation*}
d = \min_{\substack{\nb \in \ZZ^3 \setminus \{ 0 \}\\ t \in \RR^{\geq 0}} } 
              \| \qb_i + \widetilde{L}_t \nb - \qb_i \|
\geq \min_{\substack{\nb \in \ZZ^3 \setminus \{ 0 \} \\ 
\rstr \in \widehat{\Omega} }}   \| S e^{\rstr} V^{-1} \nb\|.
\end{equation*} 
Using the boundedness of
$\widehat{\Omega},$ we can limit the search to a small number of $\nb
\in \ZZ^3,$ and minimization over $\widehat{\Omega}$ leaves a quick
computation.  
For the matrices in Section~\ref{sec:algo}, if the vectors of $S$
are orthogonal, the minimum distance is
$d \approx 0.8198 a,$ where we recall that $a^3$ is the volume of
the simulation box.

\section{Numerics}

In the following, we test the consistency of our algorithm by
comparing computations for a WCA fluid under three-dimensional
elongation to those presented in~\cite{bara95}.  In previous works, the 
simulation time was restricted by the elongation of the unit cell,
though the authors in~\cite{bara95} proposed a doubling scheme 
that increased the size of the unit cell to increase the simulation
time.  This came at the tradeoff of additional computational cost.
In the following, we show that our simulations using the generalized KR boundary conditions
converge to the same macroscopic quantities even after several cell
resets.

We use the WCA potential~\cite{week71}, which is given by
\begin{equation*}
\phi(r) = \begin{cases}
\displaystyle 4 \left[ \frac{1}{r^{12}} - \frac{1}{r^6} \right] + 1, &r \leq
2^{1/6}, \\
0, &r > 2^{1/6}.
\end{cases}
\end{equation*} 
We simulate $N=512$ particles at the scaled temperature $T=0.722$ and
fluid density $\rho = 0.8442.$  For consistency with previous
works~\cite{bara95, hunt10}, we employ the SLLOD equations of
motion~\cite{evan84} with Gaussian (isokinetic)
thermostat~\cite{evan07}, which is given by
\begin{equation*}
\begin{split}
\frac{d \qb}{dt} &= \vb, \\
\frac{d \vb}{dt} &= M^{-1} \fb + A A \qb - \alpha (\vb - A \qb), \\  
\alpha &= \frac{(M^{-1} \fb - A \vb + A A \qb) \cdot (\vb - A \qb)}{(v - A \qb) \cdot (\vb - A \qb)},
\end{split}
\end{equation*} 
where $\qb \in \RR^{3N}$ denotes the vector of all particle positions,
$\vb \in \RR^{3N}$ denotes the corresponding velocities, and $\fb \in
\RR^{3N}$ denotes the interaction forces on the particles.  The factor $\alpha$
ensures that the relative kinetic energy $ \frac{1}{2} (\vb - A
\qb)^2$ is exactly preserved by the dynamics.  

\begin{figure}
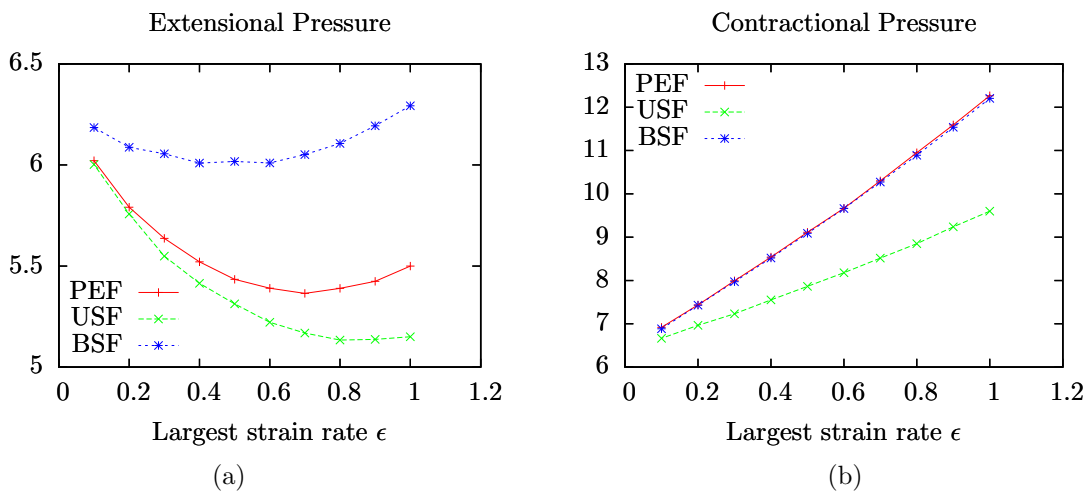

\centerline{ 
\input{figs/pressure_vs_extension.tex} 
            \input{figs/pressure_vs_extension_2.tex} }   
\centerline{ (a)\hspace{3in}  (b)} 
\caption{\label{fig:pressures}Pressures for PEF, USF, and BSF flows.
Components of the pressure tensor (which is the negative stress~\eqref{virial})
are plotted against the largest magnitude component of the velocity gradient
tensor.  In (a) the pressure in the direction of extension is plotted, while
in (b) the pressure in the direction of contraction is plotted.  These plots
show close agreement to the plots~\cite[Fig.~8 and Fig.~9]{bara95}.}
\end{figure}

We run our simulations up to time $t_{\rm max}=20,$ with time step $\Delta t =
0.002.$  The initial positions are on a lattice with random velocities
that are scaled so the system has the temperature $T=0.722.$  We allow
a decorrelation step from the initial conditions up to time $t=2,$ and
then average the desired observables until $t_{\rm max}.$  For the
largest strains, the unit cell is remapped approximately $15$
times over the course of the simulation.  We run ten realizations for
each type of flow.  We compute the virial stress tensor, 
\begin{equation}
\label{virial}
{ {\sigma}} = - \frac{1}{\det L_t} \sum_{i =1}^{N} 
\left(M (\vb_i- A \qb_i) \otimes
(\vb_i- A \qb_i) + \frac{1}{2} \sum_{\substack{i,j = 1\\j \neq
i}}^{N} ( \qb_i - \qb_j)
\otimes f^{(i j)}\right)
\end{equation}
where 
$$f^{(i j)} = - \phi'(|\qb_i - \qb_j|) \frac{\qb_i - \qb_j}{| \qb_i - \qb_j|}.$$
We also use the pressure tensor, $P = -  \sigma.$
In Figure~\ref{fig:pressures}, we plot the pressures for three
different elongational flow types, planar elongational flow (PEF), 
uniaxial stretching flow (USF), and biaxial stretching flow (BSF),
which have the respective velocity gradients
\begin{equation*}
A_{PEF} =
\left[ \begin{array}{rrr}  
\eps & & \\
& -\eps  & \\
& & 0 \\
\end{array}
\right]
\quad
A_{USF} =
\left[ \begin{array}{rrr}  
\eps & & \\
& -\eps / 2 & \\
& & -\eps / 2 \\
\end{array}
\right]
\quad 
A_{BSF} =
\left[ \begin{array}{rrr}  
- \eps & & \\
& \eps / 2 & \\
& & \eps / 2 \\
\end{array}
\right]
\end{equation*}
where $\eps > 0.$ 
In Figure~\ref{fig:pressures}(a) the pressure in the extensional direction is
plotted versus $\eps$, and in Figure~\ref{fig:pressures}(b) the pressure in the
compression direction is plotted versus $\eps$. These plots show close
agreement to the plots~\cite[Fig.~8 and Fig.~9]{bara95}.

\begin{figure}
\centerline{\input{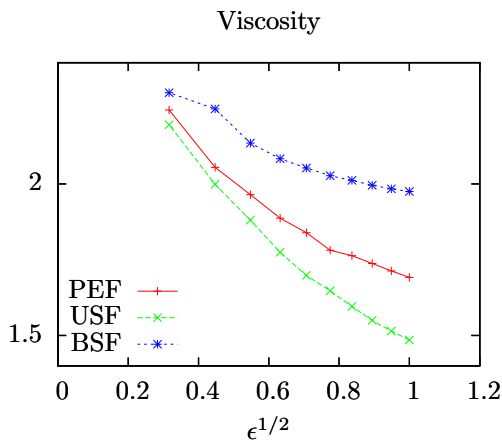}}

\caption{\label{fig:viscosity}Viscosity for PEF, USF, and BSF flows.
 These plots show close agreement
to the plots~\cite[Fig.~6]{bara95}.}
\end{figure}

For a given velocity gradient $A,$ we define $ \gamma = A + A^T,$
and define the generalized viscosity~\cite{houn92}
\begin{equation*}
 \eta = \frac{ \sigma :  \gamma}{ \gamma :  \gamma},
\end{equation*}
where $A : B = \sum_{i,j} A_{ij} B_{ij}$ denotes the 
contraction product of a pair of tensors.  In
Figure~\ref{fig:viscosity} we plot the viscosity against the square
root of $\eps.$

\begin{remark}
We note that the WCA fluid we simulate is a simple fluid, with short
decorrelation time, so that it is possible to use finite duration
simulations.  Our algorithm has more practical application for complex
molecular systems where the decorrelation time is longer than allowed
by traditional, time-restricted simulations.  The above numerics
are to show consistency of the computational results in a simple case.
\end{remark}

\section{Conclusion}

We have generalized the KR boundary conditions to handle all
homogeneous, incompressible three dimensional flows whose velocity gradient
is a nondefective matrix.  In particular, it can treat the cases of
uniaxial or biaxial flow, which could not be treated with the original
KR boundary conditions.  The boundary conditions allow the simulations
to continue for arbitrarily long times, which is important for the
simulation of complex fluids with large decorrelation times.

\section*{Acknowledgements}
The author would like to thank Gabriel Stoltz for a careful reading of
an early manuscript, as well as Bob Kohn for helpful discussions.

\end{document}